\documentstyle[11pt]{article}
\newlength{\standardunitlength}
 \newtheorem{lemma}{Lemma}
\newtheorem{theorem}{Theorem} 
\newenvironment{proof}{\noindent {\sc Proof:}}{$\Box$ \vspace{2 ex}}
\begin{document}

\begin{center}
Stein's Method and Non-Reversible Markov Chains
\end{center}

\begin{center}
By Jason Fulman
\end{center}

\begin{center}
Current affiliation: University of Pittsburgh
\end{center}

\begin{center}
Current email: fulman@math.pitt.edu
\end{center}

\begin{abstract}
  Let $W(\pi)$ be either the number of descents or inversions of a
  permutation $\pi \in S_n$. Stein's method is applied to show that
  $W$ satisfies a central limit theorem with error rate
  $n^{-1/2}$. The construction of an exchangeable pair $(W,W')$ used
  in Stein's method is non-trivial and uses a non-reversible Markov
  chain.
\end{abstract}

\section{Introduction}

        We begin by recalling two permutation statistics on the
symmetric group $S_n$ which are of interest to combinatorialists and
statisticians. A good introduction to the combinatorial aspects of
permutation statistics is Chapter 1 of Stanley \cite{Stanley}, and a superb
account of their applications to statistical problems is Chapter 6
of Diaconis \cite{Diaconisbook}.

The first statistic on $S_n$ is $Des(\pi)$, the number of descents of
$\pi$. This is defined as the number of pairs $(i,i+1)$ with $1 \leq i
\leq n-1$ such that $\pi(i) > \pi(i+1)$. Writing $\pi$ in two-line
form, this is the number of times the value of the permutation $\pi$
decreases. (A more general definition of descents exists for Coxeter
groups: the number of height one positive roots sent to negative roots
by $\pi$). The number of permutations $\pi$ in $S_n$ with $k+1$
descents is also called the Eulerian number $A(n,k)$ and has been
studied extensively \cite{DiacPiNotes}, \cite{Foata}, \cite{Knuth}.
Several proofs are known for the asymptotic ($n \rightarrow \infty$)
normality of $A(n,k)$. See for instance Diaconis and Pitman
\cite{DiacPiNotes}, Pitman \cite{PitHarper}, Bender \cite{Bender}, and
Tanny \cite{Tanny}. A proof using the method of moments should also
work.

A second well-studied statistic on $S_n$ is $Inv(\pi)$, the number of
inversions of $\pi$. In the statistics community this is called
Kendall's tau. $Inv$ is defined as the number of pairs $(i,j)$ with $1
\leq i < j \leq n$ such that $\pi(i) > \pi(j)$.  Writing $\pi$ in
two-line form, this is the number of pairs $(i,j)$ whose values are
out of order. $I(\pi)$ is also the length of $\pi$ in terms of the
standard generators $\{(i,i+1): 1 \leq i \leq n-1\}$ for $S_n$. (For
an arbitrary Coxeter group, $Inv(\pi)$ is the number of positive roots
sent to negative roots by $\pi$). Proofs of the asymptotic normality
of $Inv(\pi)$ for $S_n$ can be found in Bender \cite{Bender} and
Chapter 6 of Diaconis \cite{Diaconisbook}.

        The following definition generalizes both of these statistics. Let $M=(M_{i,j})$
be a real, anti-symmetric, $n*n$ matrix. Let $X$ be the random variable on $S_n$
defined by $X(\pi)=
\sum_{i<j} M_{\pi(i),\pi(j)}$. Setting $M_{i,j}=-1$ if $j=i+1$, $M_{i,j}=1$ if
$j=i-1$, and $M_{i,j}=0$ otherwise leads to
$X(\pi)=2Des(\pi^{-1})-(n-1)$. Setting $M_{i,j}=-1$ if $i<j$,
$M_{i,j}=+1$ if $i>j$, and $M_{i,i}=0$ leads to
$X(\pi)=2Inv(\pi^{-1})-{n \choose 2}$. Define $W=\frac{X}
{\sqrt{Var(X)}}$, so that $W$ has mean 0 and variance 1.

Charles Stein developed a method for bounding the sup
norm between the distribution of a random variable and the standard
normal distribution. 
His technique has come to be known as Stein's method. Stein's book
\cite{Ste} and the papers in this volume are good references. 

	Let us recall some notation from probability theory. If $Y,Z$ are random variables on a probability space
$(\Omega,\it{B},P)$, we let $E(Y)$ denote the expected value of $Y$
and $E^Z(Y)$ the expected value of $Y$ given $Z$, where both
expectations are taken under $P$. In the case at hand, $\Omega$ is $S_n$,
$\it{B}$ is all subsets of $S_n$, and $P$ is the uniform
distribution. Call $W,W'$ an exchangeable pair of random variables
on $S_n$ if $P(W=w_1,W'=w_2)=P(W=w_2,W'=w_1)$.

	Theorem \ref{Stein} is due to Rinott and Rotar.

\begin{theorem} \label{Stein} (\cite{RR}) Let $W,W'$ be an exchangeable pair of
real random variables such that $E^W W'=(1-\lambda)W$ with
$0<\lambda<1$. Suppose moreover that $|W'-W| \leq A$ for some constant $A$. Then for all real $x$,

\[ |P\{W \leq x\}-\Phi(x)| \leq \frac{12}{\lambda} \sqrt{Var(E^W(W'-W)^2)} + 48 \frac{A^3}{\lambda} + 8 \frac{A^2}{\sqrt{\lambda}}
\] where $\Phi$ is the standard normal distribution.

\end{theorem}

        Theorem $\ref{Stein}$ will be used to prove Theorem $\ref{Main}$.

\begin{theorem} \label{Main} Let $Des(\pi)$ and $Inv(\pi)$ be the number of
descents and inversions of $\pi \in S_n$. Then for all real $x$,

 \[ |P \{ \frac{Des-\frac{n-1}{2}}{\sqrt{\frac{n+1}{12}}} \leq x \} - \Phi(x)| \leq
\frac{C}{n^{\frac{1}{2}}} \]

 \[ |P \{ \frac{Inv-\frac{{n \choose 2}}{2}}{\sqrt{\frac{n(n-1)(2n+5)}{72}}} \leq x
\} -
\Phi(x)|
\leq
\frac{C}{n^{\frac{1}{2}}} \] where $C$ is a constant independent of $n$.

\end{theorem}

	We remark that Theorem \ref{Main} is known by other proof
techniques (see \cite{DiacPiNotes} for the case of descents and
\cite{Bic} for inversions). We recently learned that there is some
overlap with results in \cite{BCLZ}, which gives bounds for
permutation statistics using reversible Markov chains together with
Bolthausen's variation of Stein's method.
 
Section 2 shows how, for $W=Des$ or $W=Inv$, to construct an
exchangeable pair $(W,W')$ such that $E^W W'=(1-\frac{2}{n}) W$. This
step, which is usually the easy part of applying Stein's method, is
non-trivial and uses a non-reversible Markov chain equivalent to the
``move to front'' chain. The only other example in the literature in
which exchangeability was not obvious is the paper of Rinott and Rotar
\cite{RR}. A connection with this work will be mentioned in Section
2. Section 3 develops bounds for the terms on the right-hand side of
Theorem 1, and indicates why a somewhat weaker version of Theorem \ref{Stein} due to Stein can only give $n^{-1/4}$ rates.

We remark that the move to front rule on the symmetric group is a very
special case of a theory of random walk on the chambers of real
hyperplane arrangements \cite{BHR}. The corresponding Markov chains
are non-reversible and have real eigenvalues. These nonreversible
chains have recently been related to a reversible Markov chain on the
set of irreducible representations of the symmetric group
\cite{F1},\cite{F2}.

\section{Construction of an Exchangeable Pair $(W,W')$}

This section constructs $W'$ so that $(W,W')$ is an exchangeable pair
with nice properties. In most applications of Stein's method (e.g. the
examples in Stein \cite{Ste}), it is clear how to define $W'$ and
exchangeability comes for free. The situation here is more subtle.

        This being said, define $W'=W'(\pi)$ as follows. Pick $I$
uniformly at random between $1$ and $n$ and define $\pi'$ as $(I,I+1,\cdots,n)
\pi$, where $(I,I+1,\cdots,n)$ cycles by mapping $I \rightarrow I+1
\rightarrow \cdots \rightarrow n \rightarrow I$, and where permutation
multiplication is from left to right. For example, suppose that $n=7$
and $I=3$. Then the permutation $\pi$ which in 2-line form is:

\[ \begin{array}{c c c c c c c c c}
                i \  & : & 1 & 2 & 3 & 4 & 5 & 6 & 7 \\
                \pi(i) & : & 6 & 4 & 1 & 5 & 3 & 2 & 7
          \end{array}  \] is transformed to:

\[ \begin{array}{c c c c c c c c c}
                i \  & : & 1 & 2 & 3 & 4 & 5 & 6 & 7 \\
                \pi'(i) & : & 6 & 4 & 5 & 3 & 2 & 7 & 1
          \end{array}  \]

        In other words, one moves the number in position $I$ in the
second row of $\pi$ to the end of this second row. Now define
$W'(\pi)=W(\pi')$. Before discussing exchangeability, we prove Lemma
$\ref{Lambda}$, which was the motivation for the definition of $W'$
and shows that one can take $\lambda=\frac{2}{n}$ in Theorem
$\ref{Stein}$.

\begin{lemma} \label{Lambda} $E^W W'=(1-\frac{2}{n})W$
\end{lemma}

\begin{proof}
        Letting $i$ be the value of the random variable $I$, ones sees from the
definition of $W'$ that:

\begin{eqnarray*}
E^{\pi}(W'-W) & = & \frac{1}{\sqrt{Var(X)}} \frac{1}{n} \sum_{i=1}^n \sum_{j: j>i}
-2 M_{\pi(i),\pi(j)}\\ & = & \frac{1}{\sqrt{Var(X)}} \frac{1}{n} \sum_{1\leq i<j
\leq n} - 2M_{\pi(i),\pi(j)}\\ & = & -\frac{2}{n} W.
\end{eqnarray*} Since $E^{\pi}(W'-W)$ depends on $\pi$ only through $W$, the
lemma follows.
\end{proof}

        Lemma \ref{Graph} establishes a condition on $(M_{i,j})$ under which the pair
$(W,W')$ is exchangeable. This condition admittedly has limited scope, but as
will be seen, holds for the cases of descents and inversions.

\begin{lemma} \label{Graph} Given a subset $S$ of $\{1,\cdots,n\}$, for each $i
  \in S$ define $a_{i,S}=\sum_{j \in S: j>i} M_{i,j}$ and
  $b_{i,S}=\sum_{j \in S: j<i} M_{j,i}$. Suppose that for all subsets
  $S$ of $\{1,\cdots,n\}$ , there is a bijection $\Theta:S \mapsto S$
  satisfing the following conditions:

\begin{enumerate}

\item For each $i \in S$, $a_{i,S}-b_{i,S}=b_{\Theta(i),S}-a_{\Theta(i),S}.$

\item For each $i \in S$, there is a bijection $\Phi_i: S-\{i\} \mapsto S-\{
\Theta(i)
\}$ such that $M_{j,k} = M_{\Phi_i(j),\Phi_i(k)}$ for all $j,k \in S-\{i\}$.

\end{enumerate} Then $(W,W')$ is an exchangeable pair of random variables.
\end{lemma}

\begin{proof}
        It will be shown that $P\{W=a,W'=b\}=P\{W=b,W'=a\}$. For this we prove
the stronger claim that if $T=\{\pi \in S_n: \pi(j)=z_j \ for \ 1 \leq j \leq I-1
\}$, then

\[ P\{W=a,W'=b|I,\pi \in T\} =  P\{W=b,W'=a,|I,\pi \in T\} \]

        In other words, assume that the value of $I$ and the images of
$\{1,\cdots,I-1\}$ under $\pi$ are given. Let $S = \{ \pi(I), \cdots, \pi(n) \}$
be as in the hypotheses of the lemma. Now define a bijection
$\Lambda: T \mapsto T$ as follows:

\begin{enumerate}
\item $\Lambda(\pi)(j)=\pi(j)$ for $1 \leq j \leq I-1$
\item $\Lambda(\pi)(I)=\Theta(\pi(I))$
\item $\Lambda(\pi)(j)=\Phi_{\pi(I)}(\pi(j))$ for $I+1 \leq j \leq N$
\end{enumerate}

        We only show that $W(\pi)=W(\Lambda(\pi)')$, the argument that
$W(\pi')=W(\Lambda(\pi))$ being similar. Since $\pi$ and
$\Lambda(\pi)'$ agree on $1,\cdots,I-1$, it is enough to show that \[
\sum_{I<j \leq n} M_{\pi(I),\pi(j)} + \sum_{I<i<j\leq n}
M_{\pi(i),\pi(j)} = \sum_{I \leq i<j<n}
M_{\Lambda(\pi)'(i),\Lambda(\pi)'(j)} + \sum_{I \leq i<n}
M_{\Lambda(\pi)'(i),\Lambda(\pi)'(n)}. \] Now observe
that \begin{eqnarray*} \sum_{I \leq i<j<n}
M_{\Lambda(\pi)'(i),\Lambda(\pi)'(j)} & = & \sum_{I<i<j \leq n} M_{\Lambda(\pi)(i),\Lambda(\pi)(j)} \\ & = & \sum_{I<i<j \leq n}
M_{\Phi_{\pi(I)}(\pi(i)),\Phi_{\pi(I)}(\pi(j))}\\ & = & \sum_{I<i<j \leq n}
M_{\pi(i),\pi(j)}. \end{eqnarray*} The second equality is from the
definition of $\Lambda(\pi)$ and the third equality is from condition
2 in the lemma. Also observe that \begin{eqnarray*} \sum_{I \leq i< n}
M_{\Lambda(\pi)'(i),\Lambda(\pi)'(n)} & = & \sum_{I<j \leq n} M_{\Lambda(\pi)(j),\Lambda(\pi)(I)}\\ & = & \sum_{I<j \leq n}
M_{\Lambda(\pi)(j),\Theta(\pi(I))}\\ & = &
b_{\Theta(\pi(I)),S}-a_{\Theta(\pi(I)),S}\\ & = &
a_{\pi(I),S}-b_{\pi(I),S}\\ & = & \sum_{I<j \leq n}
M_{\pi(I),\pi(j)}. \end{eqnarray*} The third equality holds because $\{\Lambda(\pi)(j):I < j \leq n \} = S-\Theta(\pi(I))$. The fourth equality is from condition 1 in the lemma.
\end{proof}

{\bf Remarks}

\begin{enumerate}

\item Let us illustrate the proof of Lemma \ref{Graph} by example for $X(\pi) =
2Des(\pi^{-1})-(n-1)$. Recall that here $M_{i,j}=-1$ if $j=i+1$,
$M_{i,j}=1$ if
$j=i-1$, and $M_{i,j}=0$ otherwise. Suppose that $I=3$ and $\pi(1)=6, \pi(2)=4$.
Thus $T= \{ \pi \in S_n : \pi(1)=6, \pi(2)=4 \}$. Note that $S=\{1,2,3,5,7\}$,
because these are the images of
$\pi(j)$ for $j \geq I=3$. One observes that the bijection $\Theta:
S \mapsto S$ defined by $\Theta(1)=3$, $\Theta(2)=2$, $\Theta(3)=1$,
$\Theta(5)=5$, $\Theta(7)=7$ satisfies condition 1 of Lemma \ref{Graph} (in
general, one defines $\Theta$ by reversing within each group of consecutive
numbers in $S$). For each $i \in S$ it is also necessary to define bijections
$\Phi_i$ such that condition 2 of Lemma \ref{Graph} holds. This can be done by
pairing the elements of $S-\{i\}$ and $S-\{\Theta(i)\}$ so as to preserve their
relative order. For instance, $\Phi_1: \{2,3,5,7\} \mapsto \{1,2,5,7\}$ is
defined by
$\Phi_1(2)=1,
\Phi_1(3)=2, \Phi_1(5)=5, \Phi_1(7)=7$.

        These choices determine the bijection $\Lambda: T \rightarrow T$
constructed in Lemma \ref{Graph}. For example,

\[ \begin{array}{c c c c c c c c c c c c c c c c c c c c c c }
                i \  & : & 1 & 2 & 3 & 4 & 5 & 6 & 7 & & & & & & i \ \ \ &1 &2 &3 &4 &5 &6 &7 \\
                \pi(i) & : & 6 & 4 & 1 & 5 & 3 & 2 & 7 & & & & & & \Lambda(\pi)(i) & 6
&4&3&5&2&1&7
          \end{array}  \]

\[ \begin{array}{c c c c c c c c c c c c c c c c c c c c c c }
                i \  & : & 1 & 2 & 3 & 4 & 5 & 6 & 7 & & & & & & i \ \ \ &1 &2 &3 &4 &5 &6 &7 \\
                \pi'(i) & : & 6 & 4 & 5 & 3 & 2 & 7 & 1 & & & & & & \Lambda(\pi')(i) & 6
&4&5&2&1&7&3
          \end{array}  \]
 
        One checks that $X(\pi)=X(\Lambda(\pi'))=0$ and $X(\pi')=X(\Lambda(\pi))=2$.

\item Let us illustrate the proof of Lemma \ref{Graph} by example for $X(\pi) =
2Inv(\pi^{-1})-{n \choose 2}$. Here $M_{i,j}=-1$ if $i<j$,
$M_{i,j}=+1$ if $i>j$, and $M_{i,i}=0$. As for the case of descents, suppose
that $I=3$ and $\pi(1)=6, \pi(2)=4$. Then $T= \{ \pi \in S_n : \pi(1)=6, \pi(2)=4
\}$ and $S=\{1,2,3,5,7\}$. The bijection $\Theta: S \mapsto S$ must be defined
differently from the descent case so that condition 1 of Lemma \ref{Graph} holds.
It is easy to see that reversing the elements of $S$ works. Thus $\Theta(1)=7$,
$\Theta(2)=5$, $\Theta(3)=3$, $\Theta(5)=2$, and $\Theta(7)=1$. Defining the maps
$\Phi_i$ as in the descent case, condition 2 of Lemma \ref{Graph} holds. 

These choices determine the bijection $\Lambda: T \rightarrow T$
constructed in Lemma \ref{Graph}. For example,

\[ \begin{array}{c c c c c c c c c c c c c c c c c c c c c c }
                i \  & : & 1 & 2 & 3 & 4 & 5 & 6 & 7 & & & & & & i \ \ \ &1 &2 &3 &4 &5 &6 &7 \\
                \pi(i) & : & 6 & 4 & 1 & 5 & 3 & 2 & 7 & & & & & & \Lambda(\pi)(i) & 6
&4&7&3&2&1&5
          \end{array}  \]

\[ \begin{array}{c c c c c c c c c c c c c c c c c c c c c c }
                i \  & : & 1 & 2 & 3 & 4 & 5 & 6 & 7 & & & & & & i \ \ \ &1 &2 &3 &4 &5 &6 &7 \\
                \pi'(i) & : & 6 & 4 & 5 & 3 & 2 & 7 & 1 & & & & & & \Lambda(\pi')(i) & 6
&4&3&2&1&5&7
          \end{array}  \]

        One checks that $X(\pi)=X(\Lambda(\pi'))=1$ and $X(\pi')=X(\Lambda(\pi))=9$.

\item The above examples show that the pair $(W,W')$ is exchangeable for
descents and inversions. An interesting problem is to classify the matrices
$(M_{i,j})$ such that the pair $(W,W')$ is exchangeable. It would also be
useful to construct exchangeable pairs $(W,W')$ for other Coxeter groups.

\item Lemma 1.1 of \cite{RR} states the following. {\it Suppose that
    $\{ T^t \}$ is a stationary, nonnegative, integer valued process
    satisfying $T^{t+1}-T^{t}=+1,0$ or $-1$. Then $(T^{t},T^{t+1})$ is
    an exchangeable pair.}
  
  For the case of descents, this
  gives an alternate proof that $W,W'$ as we have defined them are an
  exchangeable pair, even though the underlying chain on permutations is not reversible. To see this, let $R^{0}$ be a uniformly
  distributed element of $S_n$; then given $R^{i}$, move to $R^{i+1}$
  according to the move random to end rule defined in the beginning of this
  section. This process is stationary. Defining $T^{t}$ to be the
  number of descents of $R^{t}$, one sees that the conditions of the
  lemma hold.

	It is interesting to note that Lemma 1.1 of \cite{RR} was
	applied there to study $W$ equal to the number of ones in a
	random pick from the stationary distribution of the antivoter
	model. The antivoter chain is not reversible, but their lemma
	implies that if $W'$ is the number of ones after a step from
	the antivoter chain, then $(W,W')$ is an exchangeable pair.

\end{enumerate}
        
\section{Bounding the Error Terms} \label{Bound}

This section bounds the error terms on the right
hand side of Theorem $\ref{Stein}$. 

        We start by computing the mean and variance of $X$ and establishing a nice
property of the pair $(W,W')$. For this it is helpful to define $A_i=\sum_{j>i}
M_{i,j}$ and $B_i=\sum_{h<i} M_{h,i}$. 

\begin{lemma} \label{meanvar} $E(X)=0$ and $Var(X)=\frac{\sum_{i<j}
(M_{i,j})^2+\sum_{i=1}^n (A_i-B_i)^2}{3}$.
\end{lemma}

\begin{proof} Observe
that the random variable $X$ on $S_n$ can be written as a sum of random variables
$X_{i,j}$ on $S_n$. Defining a random variable $X_{i,j}$ on $S_n$ by

\[ X_{i,j}(\pi) = \left\{ \begin{array}{c c}
                M_{i,j} & \mbox{if $\pi^{-1}(i) < \pi^{-1}(j)$} \\
  M_{j,i} & \mbox{if $\pi^{-1}(j) < \pi^{-1}(i)$}
        \end{array}
        \right. \] one has that: 

\[X(\pi) =  \sum_{i<j} M_{\pi(i),\pi(j)} =  \sum_{\pi^{-1}(i)<\pi^{-1}(j)} M_{i,j} =  \sum_{i<j} X_{i,j}(\pi).\]

        The mean of $X$ is $0$ since each $X_{i,j}$ has mean 0 and
expectation is linear.

        The variance of $X$ is equal to $E[(\sum_{i<j} X_{i,j}(\pi))^2]$. The
terms $E(X_{i,j}^2)$ contribute $(M_{i,j})^2$ each and thus $\sum_{i<j}
(M_{i,j})^2$ in total. The terms $E(2X_{i,j}X_{k,l})$ vanish if
$i,j,k,l$ are distinct, by independence. Now consider what
happens when two of these four indices are equal. Terms of the form
$2E(X_{i,j}X_{i,l})$ contribute $\frac{2}{3} M_{i,j}M_{i,l}$ each. The
sum of all such terms can be rewritten as $\frac{1}{3}[\sum_i
A_i^2-\sum_{i<j} (M_{i,j})^2]$. Similarly, terms of the form
$2E(X_{i,l}X_{k,l})$ contribute $\frac{1}{3}[\sum_i B_i^2-\sum_{i<j}
(M_{i,j})^2]$. Finally, terms of the form $2E(X_{i,j}X_{j,k})$ contribute
$-\frac{2}{3} M_{i,j} M_{j,k}$ each, and hence a total of
$-\frac{2}{3} \sum_i A_iB_i$. The lemma follows.
\end{proof}

As a consequence of Lemma \ref{meanvar}, one recovers the known facts that for a random permutation on n symbols, $Var(Des(\pi))=\frac{n+1}{12}$ and $Var(Inv(\pi))=\frac{n(n-1)(2n+5)}{72}$. Note that Lemma $\ref{meanvar}$ has written $Var(X)$ as a sum of positive quantities. 

\begin{lemma} $E(W'-W)^2=\frac{4}{n}$ \label{remarkable}
\end{lemma}

\begin{proof}
\begin{eqnarray*}
E(W'-W)^2 & = & E(E^W(W'-W)^2)\\ & = & E \left(E^W((W')^2+W^2-2WW')
\right)\\ & = & E \left((W')^2+E(W^2)-2 W E^W(W') \right)\\ & = & 2
Var(W) -E(2WE^W(W'))\\ & = & \frac{4}{n} Var(W)\\ & = &
\frac{4}{n}. \end{eqnarray*} The fourth equality used the fact that
$W$ and $W'$ have the same distribution. The fifth equality used Lemma
\ref{Lambda}. \end{proof}

        Lemma \ref{condit} establishes a well known inequality. For completeness, we include a proof. 

\begin{lemma} \label{condit} $E[E^W(W'-W)^2]^2 \leq E[E^{\pi}(W'-W)^2]^2$.
\end{lemma}

\begin{proof}
        Jensen's inequality says that if $g$ is a convex function, and
$Z$ a random variable, then $g(E(Z)) \leq E(g(Z))$. There is also a
conditional version of Jensen's inequality (Section 4.1 of Durrett
\cite{Durrett}) which says that if
$F$ is any $\sigma$ subalgebra of $B$, then

        \[ E(g(E(Z|F))) \leq E(g(Z)). \] The lemma follows by applying this inequality to the case $g(t)=t^2$,
$Z=E^{\pi}(W'-W)^2$, $B$ is all subsets of $S_n$, and $F$ is the $\sigma$
subalgebra of
$B$ generated by the level sets of
$W$.
\end{proof}

	Now we prove Theorem \ref{Main}.

\begin{proof} (of Theorem \ref{Main}) We will apply Theorem \ref{Stein}. Note that
 the move random to end rule changes the number of descents by at most
 one. Hence the corresponding pair $(W,W')$ satisfies $|W'-W| \leq
 \frac{2}{\sqrt{Var(X)}}$. Similarly the move random to end rule changes the
 number of inversions by at most n-1. Hence the corresponding pair
 $(W,W')$ satisfies $|W'-W| \leq \frac{2(n-1)}{\sqrt{Var(X)}}$. Thus in both
 cases $|W'-W|$ is at most $A n^{-1/2}$ for an absolute constant
 $A$. Also note by Lemma \ref{Lambda} that $E^W(W')=(1-\lambda)W$ with
 $\lambda=\frac{2}{n}$.

Thus by Theorem \ref{Stein} the result will follow if it can be shown
	that $Var(E^W(W'-W)^2) \leq \frac{B}{n^3}$. Lemma \ref{condit}
	implies that $Var(E^W(W'-W)^2) \leq Var(E^{\pi}(W'-W)^2).$
	Hence we show that $ Var(E^{\pi}(W'-W)^2) \leq \frac{B}{n^3}$.
	
Observe that
\begin{eqnarray*}
 E^{\pi}(W'-W)^2 &= &\frac{1}{Var(X)} \frac{4}{n} \sum_{i=1}^n
(\sum_{j>i} -M_{\pi(i),\pi(j)})^2\\ & = & \frac{1}{Var(X)}
\frac{4}{n} \left( \sum_{i=1}^n \sum_{j>i} (M_{\pi(i),\pi(j)})^2 + 2
\sum_{i=1}^n \sum_{i<j_1<j_2 \leq n} M_{\pi(i),\pi(j_1)}
M_{\pi(i),\pi(j_2)} \right) . \end{eqnarray*} Since $\sum_{i=1}^n
\sum_{j>i} (M_{\pi(i),\pi(j)})^2$ is independent of $\pi$, it follows
that

\begin{eqnarray*}
Var(E^{\pi}(W'-W)^2) & = & \frac{64}{Var(X)^2n^2} [\sum_{1 \leq i < j_1 <j_2 \leq
n} Var (M_{\pi(i),\pi(j_1)} M_{\pi(i),\pi(j_2)}) \\
& & +
\sum_{i < j_1 < j_2, k < l_1 <l_2 \atop (i,j_1,j_2) \neq (k,l_1,l_2)}
Cov(M_{\pi(i),\pi(j_1)} M_{\pi(i),\pi(j_2)}, M_{\pi(k),\pi(l_1)}
M_{\pi(k),\pi(l_2)})]
\end{eqnarray*}

Let us analyze this bound for the case of descents (i.e. $M_{i,j}=-1$
if $j=i+1$, $M_{i,j}=1$ if $j=i-1$, and $M_{i,j}=0$ otherwise). We
first study the summands and then divide by $Var(X)^2n^2$. The first
summand has $O(n^3)$ terms, each contributing $O(n^{-2})$; hence it is
$O(n)$. The covariance terms are also $O(n)$. To see this, first note
that the covariance vanishes if $\{i,j_1,j_2\} \cap \{k,l_1,l_2\} =
\emptyset$, so such terms can be ignored. Suppose that $i \neq
k$. Then there are $O(n^5)$ terms each contributing $O(n^{-4})$. If
$i=k$ there are subcases to consider based on which (if any) of
elements of $\{j_1,j_2\}$ are equal to elements of $\{l_1,l_2\}$. It
is straightforward to see that in all cases the contribution of the
covariance term is $O(n)$. Since $Var(X)$ is $\frac{n+1}{12}$, it
follows as desired that $ Var(E^{\pi}(W'-W)^2) \leq \frac{B}{n^3}$.

	The case of inversions is similar. The variance terms
contribute at most $O(n^3)$ and the covariance terms at most order
$O(n^5)$. Thus \[ Var(E^{\pi}(W'-W)^2) \leq \frac{B_0 n^5}{Var(X)^2n^2}
\leq \frac{B}{n^3} \] where $B_0,B$ are universal
constants. \end{proof}

	To conclude the paper, we comment on the following result of Stein \cite{Ste}.

\begin{theorem} (Stein) Let $W,W'$ be an exchangeable pair of
real random variables such that $E^W W'=(1-\lambda)W$ with
$0<\lambda<1$. Then for all real $x$,
\[ |P\{W \leq x\}-\Phi(x)| \leq 2 \sqrt{E[1-\frac{1}{2 \lambda}E^W
(W'-W)^2]^2} + (2 \pi)^{- \frac{1}{4}} \sqrt{\frac{1}{\lambda} E|W'-W|^3}
\] where $\Phi$ is the standard normal distribution.
\end{theorem} Applied to our exchangeable pair this would only yield bounds of order $n^{-1/4}$, since by Jensen's inequality $E|W'-W|^3 \geq (E(W'-W)^2)^{3/2} = \left(\frac{4}{n} \right)^{3/2}$.

\section{Acknowledgements} The author thanks to Persi Diaconis for introducing the author to Stein's method and urging him to
find applications of it to permutation statistics. We also thank
Y. Rinott for useful feedback on the 1997 version of this
paper. This research was done under the support of the National
Defense Science and Engineering Graduate Fellowship (grant no.
DAAH04-93-G-0270) and the Alfred P.  Sloan Foundation Dissertation
Fellowship.


\begin{thebibliography}{AAA}
  
\bibitem [BCLZ]{BCLZ} Bai, Z., Chao, C.C, Liang, C., and Zhao, L.,
  Error bounds in a central limit theorem of doubly indexed
  permutation matrices, {\it Annals of Stat.} {\bf 25} (1997),
  2210-2227.

  
\bibitem [Be]{Bender} Bender, E., Central and local limit theorems applied
  to asymptotic enumeration, {\it J. Combin. Theory Ser. A} {\bf 15}
  (1973), 91-111.


\bibitem [Bic]{Bic} Bickel, P., Edgeworth expansions in nonparametric statistics, {\it Annals of Stat.} {\bf 2} (1974), 1-20.
  
\bibitem [BiHaR]{BHR} Bidigare, P., Hanlon, P., and Rockmore, D., A
  combinatorial description of the spectrum for the Tsetlin library
  and its generalizations to hyperplane arrangements, {\it Duke Math
    J.} {\bf 99} (1999), 135-174.

\bibitem [Di]{Diaconisbook} Diaconis, P., {\it Group representations in
  probability and statistics}. Institute of Mathematical Statistics
  Lecture Notes, Volume 11, 1988.

\bibitem [DiP] {DiacPiNotes} Diaconis, P. and Pitman, J., Unpublished
notes on descents, 1991.

\bibitem [Du] {Durrett} Durrett, R., Probability: Theory and examples,
  Brooks/Cole Publishing Company, 1991.

\bibitem [FSc] {Foata} Foata, D. and Schutzenberger, M.P., Theorie
  geometriques des polynomes Euleriens. {\it Springer Lecture Notes in
    Math. 138} (1970).

\bibitem [F1]{F1} Fulman, J., Stein's method and Plancherel measure of the symmetric group, to appear in {\it Trans. Amer. Math. Soc.}

\bibitem [F2]{F2} Fulman, J., Card shuffling and the decomposition of tensor products, to appear in {\it Pacific J. Math}.
  
\bibitem [K] {Knuth} Knuth, D., {\it The art of computer programming.
    Volume 3. Sorting and searching.} Addison-Wesley Publishing Co.,
  1973.
  
\bibitem [P] {PitHarper} Pitman, J., Probabilistic bounds on the
  coefficients of polynomials with only real zeros. {\it J. Combin.
    Theory Ser. A} {\bf 77} (1997), 279-303.
  
\bibitem [RR] {RR} Rinott, Y., and Rotar, V., On coupling
  constructions with rates in the CLT for dependent summands with
  applications to the antivoter model and weighted $U$-statistics, {\it
  Annals Appl. Probab.} {\bf 7} (1997), 1080-1105.

\bibitem [Sta] {Stanley} Stanley, R., {\it Enumerative combinatorics, Volume
    1.}  Wadsworth and Brooks/Cole Mathematical Series, 1986.
  
\bibitem [Ste] {Ste} Stein, C., {\it Approximate computation of
    expectations.} Institute of Mathematical Statistics Lecture Notes,
  Volume 7, 1986.
  
\bibitem [T] {Tanny} Tanny, S., A probabilistic interpretation of
  Eulerian numbers. {\it Duke Math J.} {\bf 40} (1973), 717-722.

\end{thebibliography}
\end{document}